\documentclass[reqno,11pt]{amsart}  
\usepackage{amsmath} 
\usepackage{amssymb} 
\usepackage{amsthm} 
\newcommand\NoBlackBoxes{\global\overfullrule0pt} 
\NoBlackBoxes 
\theoremstyle{plain} 
 
\def\4{\kern1pt}

\def\6{\vphantom0}

\def\8{\kern-10pt} 
\def\7#1{_{(#1)}}

\makeatletter 
\let\serieslogo@\relax 
\let\@setcopyright\relax 
 
\def\speciallabelmark#1{\def\@currentlabel{#1}} 
\makeatother 
 
\begin{document} 
 
\def\ffrac#1#2{\raise.5pt\hbox{\small$\4\displaystyle\frac{\,#1\,}{\,#2\,}\4$}} 
\def\ovln#1{\,{\overline{\!#1}}} 
\def\ve{\varepsilon} 
\def\kar{\beta_r} 

\title{CONVERGENCE TO STABLE LAWS \\ 
IN RELATIVE ENTROPY 
} 
 
\author{S. G. Bobkov$^{1,4}$} 
\thanks{1) School of Mathematics, University of Minnesota, USA; 
Email: bobkov@math.umn.edu} 
\address 
{Sergey G. Bobkov \newline 
School of Mathematics, University of Minnesota  \newline  
127 Vincent Hall, 206 Church St. S.E., Minneapolis, MN 55455 USA 
\smallskip} 
\email {bobkov@math.umn.edu}  
 
\author{G. P. Chistyakov$^{2,4}$} 
\thanks{2) Faculty of Mathematics, University of Bielefeld, Germany; 
Email: chistyak@math.uni-bielefeld.de} 
\address 
{Gennadiy P. Chistyakov\newline 
Fakult\"at f\"ur Mathematik, Universit\"at Bielefeld\newline 
Postfach 100131, 33501 Bielefeld, Germany} 
\email {chistyak@math.uni-bielefeld.de} 
 
\author{F. G\"otze$^{3,4}$} 
\thanks{3) Faculty of Mathematics, University of Bielefeld, Germany; 
Email: goetze@math.uni-bielefeld.de} 
\thanks{4) Research partially supported by  
NSF grant and  
SFB 701} 
\address 
{Friedrich G\"otze\newline 
Fakult\"at f\"ur Mathematik, Universit\"at Bielefeld\newline 
Postfach 100131, 33501 Bielefeld, Germany} 
\email {goetze@mathematik.uni-bielefeld.de}


 
\begin{abstract} 
Convergence to stable laws in relative entropy is established  
for sums of i.i.d. random variables. 
\end{abstract} 

\maketitle 
\markboth{S. G. Bobkov, G. P. Chistyakov and F. G\"otze}{Entropic limit theorems}


 
 
\def\theequation{\thesection.\arabic{equation}} 
\def\E{{\bf E}} 
\def\R{{\bf R}} 
\def\C{{\bf C}} 
\def\P{{\bf P}} 
\def\H{{\rm H}} 
\def\Im{{\rm Im}} 
\def\Tr{{\rm Tr}} 
 
\def\k{{\kappa}} 
\def\M{{\cal M}} 
\def\Var{{\rm Var}} 
\def\Ent{{\rm Ent}} 
\def\O{{\rm Osc}_\mu} 
 
\def\ep{\varepsilon} 
\def\phi{\varphi} 
\def\F{{\cal F}} 
\def\L{{\cal L}} 
 
\def\be{\begin{equation}} 
\def\en{\end{equation}} 
\def\bee{\begin{eqnarray*}} 
\def\ene{\end{eqnarray*}}

 
\section{{\bf Introduction}} 
\setcounter{equation}{0} 
 
Given independent identically distributed random variables $(X_n)_{n \geq 1}$,  
consider the normalized sums 
$$ 
Z_n = \frac{X_1 + \dots + X_n}{b_n} - a_n, 
$$ 
defined for given (non-random) $a_n \in \R$ and $b_n > 0$. Assume that  
$Z_n$ converges weakly in distribution to a random variable $Z$ which has  
a non-degenerate stable law. In this paper, we would like to study 
whether or not this convergence holds in a stronger sense.  
This question has been studied and affirmatively solved in the literature,  
for example, using the total variation distance between  
the distributions of $Z_n$ and $Z$ (cf. [I-L]). 
Results in the central limit theorem suggest to consider for instance 
the stronger "entropic" distance, that is, the relative entropy, where,  
however, not so much is known. Similarly one might consider the convergence  
in terms of the closely related Fisher information (a question raised in [J], Ch.5). 
Convergence in such distances may be viewed as part of a theoretic-information 
approach to limit theorems, which has been initiated by Linnik [Li], 
who studied the behavior of the entropy of sums of 
regularized independent summands in the central limit theorem. 
 
Given random variables $X$ and $Z$ with distributions $\mu$ and $\nu$,  
respectively, 
the relative entropy of $\mu$ with respect to $\nu$, sometimes called  
informational divergence or Kullback-Leibler distance, of $\mu$ and  
$\nu$, is defined by 
$$ 
D(X||Z) = D(\mu||\nu) = \int \log \frac{d\mu}{d\nu}\,d\mu, 
$$ 
provided that $\mu$ is absolutely continuous with respect to $\nu$ 
(and otherwise $D(X||Z) = +\infty$).  
 
We consider the relative entropy with respect to so-called 
non-extremal stable laws (cf. relations (1.1) below and the definition 
before them). The aim of this note is to prove: 
 
\vskip5mm 
{\bf Theorem 1.1.} {\it Assume that the sequence of normalized sums  
$Z_n$ defined above converges weakly to a random variable $Z$ with a non-extremal 
stable limit law. Then the relative entropy  
distances converge to zero, that is, $D(Z_n||Z) \rightarrow 0$, as  
$n \rightarrow \infty$, 
if and only if $D(Z_n||Z) < +\infty$, for some $n$. 
} 
 
\vskip5mm 
In the sequel, we consider non-degenerate distributions, only. 
 
If $X_1$ has a finite second moment, weak convergence $Z_n \Rightarrow Z$ 
holds with $a_n = \sqrt{n}\ \E X_1$, $b_n = \sqrt{n}$, and $Z$ being 
normal.  
In this case, Theorem 1.1 turns into the entropic central limit theorem  
by Barron [B]; cf. also [B-C-G] for refinements and a different approach. 
Thus, Theorem 1.1 may be viewed as an extension of Barron's result. 
 
If $X_1$ has an infinite second moment, but still belongs to the domain of 
normal attraction, it follows that $D(Z_n||Z) = +\infty$ for all $n$.  
Hence, in this special case there is no convergence in relative entropy. 
 
In the remaining cases $Z$ has a stable distribution with some parameters 
$0 < \alpha < 2$ and $-1 \leq \beta \leq 1$. It then has a continuous 
density $\psi(x)$ with characteristic function $f(t) = \E\, e^{itZ}$  
described by the formula 
$$ 
\log f(t) =  
\exp\big\{iat - c|t|^\alpha \big(1 + i\beta\, {\rm sign}(t)\, 
\omega(t,\alpha)\big)\big\}, 
$$ 
where $a \in \R$, $c>0$, and $\omega(t,\alpha) = \tan(\frac{\pi \alpha}{2})$ 
in case $\alpha \neq 1$, and $\omega(t,\alpha) = \frac{2}{\pi}\,\log |t|$ 
for $\alpha = 1$. 
In particular, $|f(t)| = e^{-c|t|^\alpha}$. 
 
Exact expressions for the characteristic function are, however, not 
sufficiently informative for establishing results like Theorem 1.1. 
In order to pass from weak convergence to a stronger convergence,  
we need more information about the stable densities $\psi$. 
 
This information in turn will depend on the type of a given stable  
distribution. A stable distribution is called {\it non-extremal},  
if it is normal or, if $0<\alpha<2$  
and $-1 < \beta < 1$. In the latter case, the density $\psi$ of $Z$ is  
positive on the whole real line and satisfies asymptotic relations 
\be 
\psi(x) \sim c_0\, |x|^{-(1+\alpha)} \ \ (x \rightarrow -\infty), \qquad 
\psi(x) \sim c_1\, x^{-(1+\alpha)} \ \ (x \rightarrow +\infty) 
\en 
with some constants  $c_0,c_1>0$ (cf. [I-L], [Z]).  
 
The behavior of $\psi$ in the extremal case (when $|\beta| = 1$) is 
different. For example, when $0<\alpha<1$, the density is positive on 
a half-axis $H = (x_0,+\infty)$ or $H = (-\infty,x_0)$ of the real line, only. 
Hence, to guarantee finiteness of the relative entropies $D(Z_n||Z)$, 
one has to require that $Z_n$ take values in $H$ (which involves 
a certain requirement on the coefficients $a_n$ and $b_n$). 
Another important issue is that, as $x \rightarrow x_0$,  
$\psi(x) \rightarrow 0$ extremely fast, so the finiteness of  
$D(Z_n||Z)$ leads as well to an additional strong moment assumption  
about the distribution of $X_1$ near a point. 
A similar effect may be observed in the case $1 \leq \alpha < 2$,  
$|\beta|=1$, as well. Here $\psi$ is positive everywhere, but tends to zero  
extremely fast either near $+\infty$ or $-\infty$ 
(especially, when $\alpha = 1$). 
 
Note that the property that $X_1$ belongs to the domain of attraction  
of a stable law of index $0 < \alpha < 2$ may be expressed explicitly  
in terms of the distribution function $F(x) = \P\{X_1 \leq x\}$.  
Namely, we have $Z_n \Rightarrow Z$  
with some $b_n > 0$ and $a_n \in \R$, if and only if 
\bee 
F(x) & = & 
(c_0 + o(1))\,|x|^{-\alpha} B(|x|) \qquad (x \rightarrow -\infty), \\ 
1-F(x) & = & 
(c_1 + o(1))\,x^{-\alpha}\ B(x) \qquad \ \ (x \rightarrow +\infty), 
\ene 
for some constants $c_0,c_1 \geq 0$ that are not both zero, and 
where $B(x)$ is a slowly varying function in the sense of Karamata  
(cf. [I-L], [Z]). 
 
Furthermore, in the non-extremal case the condition $D(Z_n||Z) < +\infty$ 
in Theorem 1.1 is equivalent to saying that 
$Z_n$ has a density $p_n$ with finite entropy 
$$ 
h(Z_n) = -\int_{-\infty}^{+\infty} p_n(x) \log p_n(x)\,dx 
$$ 
(with an additional requirement that $\E X_1^2 < +\infty$, when $Z$ is normal). 
Once this property is fulfilled for a particular value 
$n = n_0$, it continues to hold for all $n \geq n_0$. 
This will be explained in the next section. 
 
We shall turn to the proof which we divide into several steps. 
For completeness, the argument will cover the normal case as well.

 
\vskip5mm 
\section{{\bf Remarks on Relative Entropy}} 
\setcounter{equation}{0} 
 
First, let us give necessary and sufficient conditions for the property 
$D(Z_n||Z) < +\infty$, especially when $Z$ has a non-extremal 
stable distribution. 
 
Given a random variable $X$ with density $p$, consider the entropy functional 
$$ 
h(X) = -\int p(x) \log p(x)\,dx. 
$$ 
In general, it may or may not be defined as a Lebesgue integral. 
(Here and below, we often omit the limits of integration when integrating  
over the whole real line.) 
 
It is well-known that, if $X$ has a finite second moment, then $h(X)$  
is well-defined, and one has an upper estimate 
\be 
h(X) \leq h(Z), 
\en 
where $Z$ is a normal random variable with the same mean and variance as $X$. 
(Here the value $h(X) = -\infty$ is possible.) This important observation 
may be generalized with respect to other (not necessarily normal) 
reference measures. 
 
\vskip5mm 
{\bf Proposition 2.1.} {\it Let $X$ denote a random variable with density $p$. 
Assume that $\psi$ denotes a probability density on the real line, such that  
$\psi(x) = 0$ implies $p(x) = 0$ a.e., and such that  
$\E \log^+(\frac{1}{\psi(X)}) < +\infty$. Then $h(X)$ exists and satisfies 
\be 
h(X) \leq \E \log\frac{1}{\psi(X)}. 
\en 
} 
 
\vskip2mm 
This bound seems to be folklor knowledge and is based on a direct application 
of Jensen's inequality. To recall the argument, let $Z$ have density $\psi$, 
and assume that $\psi(x) = 0$ implies $p(x) = 0$ a.e. The relative entropy 
given by 
$$ 
D(X||Z) = \int p(x) \log \frac{p(x)}{\psi(x)}\,dx = \E\, \xi \log \xi, 
\qquad \xi = \frac{p}{\psi}, 
$$ 
is then well-defined, where the expectation refers to the probability space 
$(\R,\psi(x)\,dx)$. Moreover, since $\E \xi = 1$, and due to the convexity  
of the function $t \rightarrow t \log t$, the expectation $\E\, \xi \log \xi$ 
exists and is non-negative. But   
$p \log \frac{p}{\psi} + p\log\frac{1}{p} = p\log\frac{1}{\psi}$ and, 
by the assumption, $\int  p\log\frac{1}{\psi}\,dx$ exists and 
does not take the value $+\infty$. Hence, 
$h(X) = \int p\log\frac{1}{p}\,dx$ also exists and 
cannot take the value $+\infty$. In addition, the equality 
\be 
D(X||Z) = -h(X) + \E \log\frac{1}{\psi(X)} 
\en 
is justified. Here, the left-hand side is non-negative, so the inequality 
(2.2) immediately follows. 
 
\vskip2mm 
Choosing for $\psi$ a normal density with the same mean and the same variance 
as for $X$, one may easily see that (2.2) reduces to (2.1). 
 
As another example, choosing for $\psi$ the density of the Cauchy measure, 
we obtain a weaker moment condition $\E \log(1 + |X|) < +\infty$,  
which guarantees that the entropy exists and satisfies $h(X) < +\infty$. 
 
Let us return to the definition of the relative entropy, 
$$ 
D(X||Z) = \int p(x) \log \frac{p(x)}{\psi(x)}\,dx, 
$$ 
assuming that $Z$ has density $\psi$. In order to describe when this  
distance is finite, one may complement Proposition 2.1 with the following: 
 
\vskip5mm 
{\bf Proposition 2.2.} {\it Let $X$ be a random variable with density $p$. 
Assume that 
 
\vskip2mm 
$a)$ $\psi(x) = 0 \Rightarrow p(x) = 0$ a.e.;  
 
\vskip0.7mm 
$b)$ $\E \log^+(\frac{1}{\psi(X)}) < +\infty$; 
 
\vskip0.7mm 
$c)$ $h(X)$ is finite. 
 
\vskip2mm 
\noindent 
Then $D(X||Z)$ is finite and is given by $(2.3)$. 
Conversely, if $\int \psi(x)^{\gamma}\,dx < +\infty$,  
for some $0 < \gamma < 1$, then the conditions $a)-c)$ are also necessary  
for the relative entropy to be finite. 
} 
 
\vskip5mm 
{\bf Proof.} Condition $a)$ means that the distribution of $X$ is absolutely 
continuous with respect to the distribution of $Z$ (which is necessary for the 
finiteness of the relative entropy). 
 
Assuming condition $a)$, we have $D(X||Z) < + \infty$, if and only if the  
integral 
$$ 
I = \int p(x) \log\bigg(1 + \frac{p(x)}{\psi(x)}\bigg) dx 
$$ 
is finite. Let $I$ indeed be finite. Split the real line into the two sets  
$A = \{x \in \R: p(x) \geq \psi(x)^{\gamma'}\}$ and $B = \R \setminus A$, 
where $\gamma < \gamma' < 1$. Restricting the integration to $A$, we get 
$$ 
I \, \geq \, \int_A p(x) \log\bigg(1 + \frac{1}{\psi(x)^{1 - \gamma'}}\bigg) dx 
  \, \geq \, (1 - \gamma')\,\int_A p(x)\,\log^+\!\big(\frac{1}{\psi(x)}\big)\,dx. 
$$ 
Hence, $\int_A\, p(x)\,\log^+(\frac{1}{\psi(x)})\,dx < +\infty$. On the other hand, 
using $t^{\gamma'}\log(1/t) \leq Ct^\gamma$ ($0 \leq t \leq 1$), we get 
\bee 
\int_B p(x)\,\log^+\!\big(\frac{1}{\psi(x)}\big)\,dx 
 & \leq & 
\int_B \psi(x)^{\gamma'} \log^+\!\big(\frac{1}{\psi(x)}\big)\,dx \\ 
 & & \hskip-10mm = \ 
\int_B \psi(x)^{\gamma'} \log\big(\frac{1}{\psi(x)}\big)\,1_{\{\psi(x) \leq 1\}}\,dx 
 \, \leq \, 
C \int_{-\infty}^{+\infty} \psi(x)^{\gamma}\,\,dx \, < \, +\infty. 
\ene 
As a result, $\int p(x)\,\log^+\!\big(\frac{1}{\psi(x)}\big)\,dx < +\infty$,  
that is, $b)$ is fulfilled. 
 
But then, by Proposition 2.1, the entropy is well-defined in the Lebesgue sense and,  
moreover, $h(X) < +\infty$. This justifies writing (2.3), which implies that  
$h(X) > -\infty$ as well. Hence, the property $c)$ holds as well. 
 
In the other direction, assuming that $a)-c)$ are fulfilled, one may start with  
equality (2.3), which shows that $D(X||Z)$ is finite.  
 
Thus, Proposition 2.2 is proved. 
 
\vskip5mm 
When $\psi$ is the density of a stable law, the condition  
$\int \psi(x)^{\gamma}\,dx < +\infty$ (for some $0 < \gamma < 1$) is fulfilled,  
so the properties $a)-c)$ are necessary and sufficient for the finiteness of  
the relative entropy with respect to $\psi$. In fact, a more detailed  
conclusion may be stated according to the types of stable laws. 
 
\vskip5mm 
{\bf Corollary 2.3.} {\it If $Z$ is normal, then $D(X||Z) < + \infty$, 
if and only if $X$ has a finite second moment and finite entropy. 
} 
 
\vskip5mm 
{\bf Corollary 2.4.} {\it If $Z$ has a non-extremal stable distribution,  
which is not normal, then $D(X||Z) < + \infty$, if and only if $X$ has a finite  
logarithmic moment $\E \log(1 + |X|)$ and finite entropy. 
} 
 
\vskip5mm 
This follows from Proposition 2.2 and the property (1.1). 
Let us recall that the condition $\E \log(1 + |X|) < +\infty$ insures that  
the entropy of $X$ exists and, moreover, $h(X) < +\infty$. 
 
The situation where $Z$ has an extremal stable distribution is a bit more  
delicate, but may be studied on the basis of Proposition~2.2 as well. 
However, we do not discuss this case here. 
 
These characteristions may be simplified for normalized sums  
$Z_n = \frac{X_1 + \dots + X_n}{b_n} - a_n$ with i.i.d. summands as in  
Theorem 1.1, provided that the sequence $Z_n$ is weakly convergent in  
distribution. Indeed, the property $Z_n \Rightarrow Z$, where $Z$ has  
a stable distribution with parameter $0 < \alpha < 2$ implies that 
$\E\, |X_1|^s < +\infty$ for any $0 < s < \alpha$, cf. [I-L], [Z].  
Hence, $\E\, |Z_n|^s < +\infty$ for all $n \geq 1$, and thus  
the random variables $Z_n$ have finite logarithmic moments. 
 
\vskip5mm 
{\bf Corollary 2.5.} {\it Assume that $Z_n$ converges weakly to a random 
variable $Z$ with a non-extremal stable limit law, which is not normal.  
Then, for each $n \geq 1$, the finiteness of the relative entropy 
$D(Z_n||Z)$ is equivalent to the finiteness of the entropy of $Z_n$. 
} 
 
\vskip5mm 
A similar conclusion holds in the normal case as well, provided that 
$b_n \sim \sqrt{n}$. Here it is well-known that $Z_n \Rightarrow Z$ 
implies that $\E\, |X_1|^2 < +\infty$. 
 
Finally, let us mention another property of non-extremal 
stable distributions. 
 
\vskip5mm 
{\bf Corollary 2.6.} {\it Assume that $Z$ has a non-extremal 
stable distribution. If the relative entropy $D(Z_n||Z)$ is finite for some  
$n=n_0$, it will be finite for all $n \geq n_0$. 
} 
 
\vskip5mm 
{\bf Proof.} 
By Jensen's inequality, $h(X+Y) \geq h(X)$ for all independent 
summands such that $h(X)$ exists and $h(X)>-\infty$. 
If $D(Z_{n_0}||Z) < +\infty$, then $h(Z_{n_0})$ is finite according to  
Proposition 2.2. Hence, for the sums $S_n = X_1 + \dots + X_n$,  
$h(S_n)$ exists for all $n \geq n_0$ and 
$$ 
h(S_n) \geq h(S_{n_0}) = \log b_{n_0} + h(Z_{n_0}) > -\infty. 
$$ 
 
Thus, $h(Z_n)$ exists with $h(Z_n) > -\infty$. 
In addition, by Corollary 2.4 (if $Z$ is not normal), we also have that  
$\E \log(1+|Z_{n_0}|) < +\infty$. By convexity of the function 
$u \rightarrow \log(1+u)$ ($u \geq 0$), this yields  
$\E \log(1+|X_1|) < +\infty$. In turn, since 
$\log(1+u+v) \leq \log(1+u) + \log(1+v)$ ($u,v \geq 0$), we get 
$\E \log(1+|S_n|) < +\infty$, for all $n$. 
In particular, $h(S_n) < +\infty$, according to Proposition 2.1.  
Therefore, $h(Z_n)$ is finite, and applying Proposition 2.2, 
we conclude that $D(Z_n||Z)$ is finite for all $n \geq n_0$. 
 
Using Corollary 2.3, a similar argument applies to the normal case as well. 
 
\vskip2mm 
Note that Corollary 2.6 does not extend to the class of extremal  
stable distributions.


\vskip5mm 
\section{{\bf Binomial Decomposition of Convolutions}} 
\setcounter{equation}{0}

\vskip2mm 
Given independent identically distributed random variables $(X_n)_{n \geq 1}$ 
and numbers $a_n \in \R$, $b_n > 0$, consider the sums  
$$ 
S_n = X_1 + \dots + X_n \quad {\rm and} \quad  
Z_n = \frac{S_n}{b_n} - a_n. 
$$ 
 
If $Z$ is a random variable with an absolutely continuous distribution 
(not necessarily stable), the condition $D(Z_n||Z) < +\infty$ ($n \geq n_0$)  
used in Theorem 1.1 implies that, for any such $n$, $Z_n$ has an  
absolutely continuous distributions with density, say $p_n(x)$.  
For simplicity, we may and will assume that $n_0=1$, that is, already $X_1$  
has a density $p(x)$. (The case, where $Z_n$ have densities starting from  
$n \geq n_0$ with $n_0 > 1$ requires minor modifications only). 
 
Since it is advantageous in the following to work with bounded densities,  
we slightly modify $p_n$ at the expense of a small change in the relative  
entropy.  
For a given number $0 < b < \frac{1}{2}$, split $H$ into the two Borel sets 
$H_1$ and $H_0$, such that $p$ is bounded by a constant $M$ on $H_1$ with 
$$ 
b = \int_{H_0} p(x)\,dx. 
$$ 
Consider the decomposition 
\be 
p(x) = (1-b) \rho_1(x) + b \rho_0(x), 
\en 
where $\rho_1$, $\rho_0$ are the normalized restrictions of $p$  
to the sets $H_1$ and $H_0$, respectively. Hence, 
for the convolutions we have a binomial decomposition 
$$ 
p^{*n} = \sum_{k=0}^n C_n^k\, (1-b)^k\, b^{n-k}\, \rho_1^{*k} * \rho_0^{*(n-k)}. 
$$ 
This function represents the density of $S_n$. 
 
For $n \geq 2$, we split the above sum into the two parts, so that 
$p^{*n} = \rho_{n1} + \rho_{n0}$ with 
$$ 
\rho_{n1} =  
\sum_{k = 2}^n C_n^k\, (1-b)^k\, b^{n-k}\, \rho_1^{*k} * \rho_0^{*(n-k)}, \qquad 
\rho_{n0} =  
b^n \rho_0^{*n} + n\, (1-b) b^{n-1}\, \rho_1 * \rho_0^{*(n-1)}. 
$$ 
Note that 
\be 
\ep_n \equiv \int \rho_{n0}(x)\,dx = b^n + n\, (1-b) b^{n-1} < nb^{n-1}. 
\en 
Finally define  
\be 
\widetilde p_n(x) = \frac{b_n}{1 - \ep_n}\,\rho_{n1}\big(a_n + b_n x\big), \qquad 
p_{n0}(x) = \frac{b_n}{\ep_n}\,\rho_{n0}(a_n + b_n x).  
\en 
Thus, for the densities $p_n$ of $Z_n$ we have the decomposition 
\be 
p_n(x) = (1 - \ep_n) \widetilde p_n(x) + \ep_n p_{n0}(x). 
\en 
 
\vskip2mm 
The (probability) densities $\widetilde p_n$ are bounded and provide a strong  
approximation for $p_n$ regardless of the choice of numbers $a_n$ and $b_n$ from  
the definition of $Z_n$. In particular, from (3.2) and (3.4) and using  
$b < \frac{1}{2}$, it follows that 
\be 
\int |\widetilde p_n(x) - p_n(x)|\,dx < 2^{-n},  
\en 
for all $n$ large enough. One of the immediate consequences of this estimate 
is the bound  
\be 
|\widetilde f_n(t) - f_n(t)| < 2^{-n} \qquad (t \in \R) 
\en 
for the corresponding characteristic functions  
$$ 
\widetilde f_n(t) = \int\, e^{itx} \widetilde p_n(x)\,dx, \qquad 
f_n(t) = \int\, e^{itx} p_n(x)\,dx. 
$$ 
 
Under mild conditions on $a_n$ and $b_n$, the approximation (3.5)  
may be sharpened by using a polynomial weight function in the $L^1$-distance.

\vskip5mm 
{\bf Lemma 3.1.} {\it If $\E\,|X_1|^s < +\infty$ $(s>0)$, and  
$|a_n| + 1/b_n = O(n^\gamma)$ with some $\gamma>0$, then for all $n$  
large enough, 
$$ 
\int |x|^s\, |\widetilde p_n(x) - p_n(x)|\,dx < 2^{-n}. 
$$ 
} 
 
{\bf Proof.} We refine arguments from the proof of a similar Lemma 2.1 in [B-C-G]. 
By (3.4),  
$$ 
|\widetilde p_n(x) - p_n(x)| \leq \ep_n (\widetilde p_n(x) + p_{n0}(x)),  
$$ 
so 
\begin{eqnarray} 
\int |x|^s\, |\widetilde p_n(x) - p_n(x)|\,dx 
 & \leq & 
\frac{\ep_n}{1 - \ep_n} \, b_n^{-s} \int |x-a_n|^s\, \rho_{n1}(x)\,dx 
 \nonumber \\ 
 & & \hskip7mm + \ b_n^{-s} \int |x-a_n|^s\, \rho_{n0}(x)\,dx. 
\end{eqnarray} 
 
Let $U_1,U_2,\dots$ be independent copies of $U$ and $V_1,V_2,\dots$  
be independent copies of $V$ (that are also independent of all $U_n$), where 
$U$ and $V$ are random variables with densities $\rho_1$ and $\rho_0$,  
respectively. From (3.2) 
\be 
\beta_s \equiv \E\, |X_1|^s = (1-b)\, \E\, |U|^s + b\, \E\, |V|^s, 
\en 
so $\E\, |U|^s \leq \beta_s/b$ and $\E\, |V|^s \leq \beta_s/b$  
(using $b < \frac{1}{2}$). Consider the sums  
$$ 
S_{k,n} = U_1 + \dots + U_k + V_1 + \dots + V_{n-k}, \qquad 
0 \leq k \leq n. 
$$ 
If $s \geq 1$, then by the triangle inequality in the space $L^s$ 
with norm $\|\xi\|_s = (\E |\xi|^s)^{1/s}$, we get 
\be 
\|S_{k,n}\|_s \ \leq \ \sum_{j=1}^k \|U_j\|_s \, +  \sum_{j=k+1}^n \|V_j\|_s 
\, = \, k \|U\|_s + (n-k) \|V\|_s \, \leq \, n\, (\beta_s/b)^{1/s}. 
\en 
Hence, $\E\, |S_{k,n}|^s \leq \frac{\beta_s}{b}\, n^s$ and,  
by Jensen's inequality, 
$\E\, |S_{k,n} - a_n|^s \leq 2^s (\frac{\beta_s}{b}\, n^s + |a_n|^s)$. 
 
If $0<s<1$, one can just use  
$$ 
|S_{k,n}|^s \leq |U_1|^s + \dots + |U_k|^s + |V_1|^s + \dots + |V_{n-k}|^s, 
$$ 
implying that $\E\, |S_{k,n}|^s \leq n \beta_s/b$ and  
$\E\, |S_{k,n} - a_n|^s \leq n\,\frac{\beta_s}{b} + |a_n|^s$. 
In both cases,  
$$ 
\E\, |S_{k,n} - a_n|^s \leq 2^s  
\bigg(\frac{\beta_s}{b}\, n^{\max(s,1)} + |a_n|^s\bigg). 
$$ 
Hence,  
\bee 
\int |x-a_n|^s\, \rho_{n1}(x)\,dx  
 & = & 
\sum_{k=2}^{n} C_n^k\, (1-b)^k\, b^{n-k}\, \E\, |S_{k,n} - a_n|^s \\ 
 & \leq &  
2^s\bigg(\frac{\beta_s}{b}\,n^{\max(s,1)} + |a_n|^s\bigg)\,(1-\ep_n), 
\ene 
\bee 
\int |x-a_n|^s\, \rho_{n0}(x)\,dx  
 & = &  
\sum_{k=0}^{1} \ C_n^k\, (1-b)^k\, b^{n-k}\, \E\, |S_{k,n} - a_n|^s \\  
 & \leq &  
2^s\bigg(\frac{\beta_s}{b}\,n^{\max(s,1)} + |a_n|^s\bigg)\, \ep_n. 
\ene 
The two estimates may be used in (3.7), and we get 
$$ 
\int |x|^s\, |\widetilde p_n(x) - p_n(x)|\,dx \, \leq \, 
\frac{2^{s+1}}{b_n^s}\, 
\bigg(\frac{\beta_s}{b}\,n^{\max(s,1)} + |a_n|^s\bigg)\, \ep_n. 
$$ 
It remains to apply (3.2) together with the assumption on $(a_n,b_n)$.  
Lemma 3.1 is proved.

\vskip5mm 
{\bf Lemma 3.2.} {\it For any $t_0>0$, there are positive constants $c$ and $C$ 
such that, for all $n \geq 2$, 
$$ 
\int_{|t| \geq t_0 b_n} |\widetilde f_n(t)|\,dt \, < \, Cb_n\, e^{-cn}. 
$$ 
} 
 
{\bf Proof.} Consider the densities 
$\rho = \rho_1^{*k} * \rho_0^{*(n-k)}$ appearing in the definition of $\rho_{n1}$. 
Their Fourier transforms (i.e., the corresponding characteristic functions)  
are connected by 
$$ 
\hat \rho(t) = \hat\rho_1(t)^k \, \hat \rho_0(t)^{n-k} \qquad (t \in \R). 
$$ 
 
By the Riemann-Lebesgue theorem, $|\hat\rho_j(t)| \leq e^{-c}$,  
for all $|t| \geq t_0$ with some constant $c>0$ ($j=0,1$). 
Hence, whenever $2 \leq k \leq n$. 
$$ 
|\hat \rho(t)| \leq A\,|\hat\rho_1(t)|^2, \qquad A = e^{-c(n-2)}. 
$$ 
By the decomposition construction, $\rho_1(x) \leq M$, for all $x$.  
Applying Plancherel's formula, we get 
$$ 
\int_{|t| \geq t_0} |\hat\rho(t)|\,dt < A \int |\hat\rho_1(t)|^2\,dt 
= 2\pi A\, \int \rho_1(x)^2\,dx \leq 2\pi AM. 
$$ 
As a consequence, the density $\rho_{n1}$ satisfies a similar inequality 
$$ 
\int_{|t| \geq t_0} |\hat\rho_{n1}(t)|\,dt < 2\pi AM\, (1 - \ep_n). 
$$ 
But, by (3.3), $\widetilde f_n(t) =  
\frac{1}{1 - \ep_n}\, e^{-ita_n/b_n} \hat\rho_{n1}(t/b_n)$, so 
$$ 
\int_{|t| \geq t_0 b_n} |\widetilde f_n(t)|\,dt =  
\frac{b_n}{1 - \ep_n} \int_{|t| \geq t_0} |\hat\rho_{n1}(t)|\,dt <  
2\pi M b_n\,e^{-c(n-2)}. 
$$ 
Thus, Lemma 3.2 is proved. 
 
\vskip5mm 
{\bf Remark 3.3.} If $Z_n \Rightarrow Z$, where $Z$ has a stable distribution 
of index $0 < \alpha \leq 2$, then necessarily  
$b_n \sim n^{1/\alpha} B(n)$, as $n \rightarrow \infty$, where $B$  
is a slowly varying function in the sense of Karamata (cf. [I-L]). 
Using standard arguments (cf. e.g. [La]), one can show as well that 
$a_n = o(n)$. Thus, the conditions of Lemma 3.1 for the 
coefficients $(a_n,b_n)$ are fulfilled, once there is a weak convergence.

 
\vskip5mm 
\section{{\bf Entropic Approximation of $p_n$ by $\widetilde p_n$}} 
\setcounter{equation}{0}

\vskip2mm 
We need to extend the assertion of Lemma 3.1 to the relative entropies with  
respect to the stable laws. Thus, assume that $Z$ has a stable distribution  
of index $\alpha \in (0,2]$ with density $\psi(x)$  
(thus including the normal law).  
Put 
$$ 
D_n = D(Z_n||Z) = \int p_n(x) \log \frac{p_n(x)}{\psi(x)}\ dx, \qquad 
\widetilde D_n = \int \widetilde p_n(x) \log \frac{\widetilde p_n(x)}{\psi(x)}\ dx, 
$$ 
where $\widetilde p_n$ are defined according to the decomposition (3.4) for  
the densities $p_n$ of the normalized sums 
$$ 
Z_n = \frac{1}{b_n}\,(X_1 + \dots + X_n) - a_n \qquad (a_n \in \R, \ b_n > 0). 
$$  
Here, as before, $X_k$ denote independent identically distributed random  
variables. 
 
In the lemma below, it does not matter whether or not the sequence 
$Z_n$ converges weakly to $Z$.

\vskip5mm 
{\bf Lemma 4.1.} {\it Assume that the distribution of $Z$ is non-extremal. 
If $D_n$ is finite for all $n \geq n_0$, and  
$|a_n| + \log b_n + 1/b_n = O(n^\gamma)$ with some $\gamma > 0$, then 
$$ 
|\widetilde D_n - D_n| < 2^{-n}, 
$$ 
for all $n$ large enough. 
} 
 
\vskip5mm 
{\bf Proof.} To simplify the notations (and the argument), we assume that  
$n_0 = 1$. In particular, $D_1 = D(X_1||Z)$ is finite, hence the entropy  
$h(X_1)$ is finite and $\E \log\frac{1}{\psi(X_1)} < +\infty$, according  
to Proposition 2.2. Define 
$$ 
D_{n0} = \int p_{n0}(x) \log \frac{p_{n0}(x)}{\psi(x)}\ dx. 
$$ 
 
By convexity of the function $L(u) = u\log u$ ($u \geq 0$)  
it follows that $D_n \leq (1 - \ep_n)\widetilde D_n + \ep_n D_{n0}$ and,  
as easy to see, 
$$ 
D_n \geq \big((1 - \ep_n)\widetilde D_n + \ep_n D_{n0}\big) + 
\ep_n \log \ep_n + (1-\ep_n) \log(1-\ep_n). 
$$ 
By the definition (3.2) of $\ep_n$, the two estimates give 
\be 
|\widetilde D_n - D_n| < Cn\, (n + \widetilde D_n + D_{n0})\, b^{n-1}, 
\en 
which holds for all $n \geq 1$ with some constant $C$. In addition, using an 
elementary inequality 
$L((1 - b)\,u + b v) \geq (1-b)\, L(u) - \frac{1}{e}\,u - \frac{1}{e}$ 
($u,v \geq 0$, $0 \leq b \leq 1$), we get from (3.1) that 
$$ 
D(X_1||Z) = \int L\bigg(\frac{p(x)}{\psi(x)}\bigg)\,\psi(x)\,dx  
\geq (1-b) \int \rho_1(x) \log \frac{\rho_1(x)}{\psi(x)}\ dx -  
\frac{2}{e}. 
$$ 
A similar inequality also holds for $\rho_0$ with $b$ in place of $1-b$, so 
$$ 
D(X_1||Z) \geq (1-b) D(U||Z) - \frac{2}{e}, \qquad  
D(X_1||Z) \geq b D(V||Z) - \frac{2}{e}, 
$$ 
where $U$ and $V$ have densities $\rho_1$ and $\rho_0$, respectively. 
Hence, these random variables have finite entropies, and by Proposition 2.2, 
\be 
\E\, \log \frac{1}{\psi(U)} < +\infty, \qquad 
\E\, \log \frac{1}{\psi(V)} < +\infty. 
\en 
 
Let $U_1,U_2,\dots$ be independent copies of $U$ and let $V_1,V_2,\dots$  
be independent copies of $V$ (which are independent of all $U_n$ as well).  
Again, by convexity of the function $u\log u$, 
\begin{eqnarray} 
\widetilde D_n 
 & \leq &  
\frac{1}{1-\ep_n} \sum_{k=2}^{n} C_n^k\, (1-b)^k\, b^{n-k} 
\int r_{k,n}(x)\, \log \frac{r_{k,n}(x)}{\psi(x)}\,dx, \\ 
D_{n0} 
 & \leq & 
\frac{1}{\ep_n} \, \sum_{k=0}^{1} \ C_n^k\, (1-b)^k\, b^{n-k} 
\int r_{k,n}(x)\, \log \frac{r_{k,n}(x)}{\psi(x)}\,dx, 
\end{eqnarray} 
where $r_{k,n}$ are the densities of the normalized sums 
$$ 
R_{k,n} = \frac{S_{k,n}}{b_n} - a_n = 
\frac{1}{b_n}\, (U_1 + \dots + U_k + V_1 + \dots + V_{n-k}) - a_n,  
\qquad 0 \leq k \leq n. 
$$ 
 
Now, write 
\be 
D(R_{k,n}||Z) = \int r_{k,n}(x)\, \log \frac{r_{k,n}(x)}{\psi(x)}\,dx 
= -h(R_{k,n}) + \int r_{k,n}(x)\, \log \frac{1}{\psi(x)}\,dx, 
\en 
using the entropy functional $h(R) = -\int r(x)\log r(x)\,dx$. 
Adding independent summands to $R$ will only increase the value of this 
functional. Hence, for any $1 \leq k \leq n$, 
\bee 
h(R_{k,n}) & = & 
-\log b_n + h(U_1 + \dots + U_k + V_1 + \dots + V_{n-k}) \\ 
 & \geq &  
-\log b_n + h(U). 
\ene 
For $k=0$ there are similar relations (with $V$ replacing $U$), so whenever 
$0 \leq k \leq n$, we have $h(R_{k,n}) \geq -\log b_n - C$ 
with some constant $C$. Inserting in (4.5), we arrive at 
\be 
D(R_{k,n}||Z) \leq \log b_n + C + \E \log\frac{1}{\psi(R_{k,n})}. 
\en 
 
\vskip2mm 
{\bf Case 1}: $Z \sim N(a,\sigma^2)$ with some $a \in \R$ and $\sigma > 0$.  
Using (4.6) and (3.8)-(3.9) with $s=2$ and the assumption $\E X_1^2 < +\infty$  
(due to the assumption $D_1 < +\infty$, cf. Corollary 2.3), we get 
\bee 
D(R_{k,n}||Z)  
& \leq &  
\log b_n + C_1 + C_2\, \E\, |R_{k,n}|^2 \\ 
 & \leq &  
\log b_n + C_3 + C_4\, \frac{n^2}{b_n^2} + C_5\, a_n^2 
\ene 
with some constants $C_j$ depending on $a,\sigma,b$ and $\E X_1^2 $.  
Using the condition on $a_n$ and $b_n$, we conclude that 
$D(R_{k,n}||Z) \leq n^{\gamma'}$ with some $\gamma'$ for all $n$  
large enough. Applying this in (4.3)-(4.4), (4.1) yields  
$|\widetilde D_n - D_n| = o(b_1^n)$, whenever $b < b_1 < \frac{1}{2}$. 
 
\vskip2mm 
{\bf Case 2}: $Z$ has a non-extremal stable distribution.  
By (4.2), we have $\E \log(1 + |U|) <+\infty$ and  
similarly for $V$. In addition, by (1.1), 
\be 
\E \log\frac{1}{\psi(R_{k,n})} \leq A + B\, \E \log(1 + |R_{k,n}|) 
\en 
with some constants $A$ and $B$. To bound the last expectation, one may use  
the inequality $\log(1 + |u_1 + u_2|) \leq \log(1 + |u_1|) + \log(1 + |u_2|)$,  
valid for all real numbers $u_1$, $u_2$, together with  
$\log(1 + u) \leq 1/u + \log u$ ($u > 0$). They yield 
$$ 
\log(1 + |R_{k,n}|) \leq \log(1 + |a_n|) + \log(1 + |S_{k,n}|/b_n), 
$$ 
while 
$\log(1 + |S_{k,n}|/b_n) \leq 1/b_n + \log(1 + |S_{k,n}|)$. From this, we get 
$$ 
\log(1 + |R_{k,n}|) \leq 1 + \log(1 + |a_n|) +  
\sum_{j=1}^k \log(1 + |U_j|) + \sum_{j=k+1}^n \log(1 + |V_j|), 
$$ 
and therefore 
\bee 
\E \log(1 + |R_{k,n}|) 
 & \leq &  
\frac{1}{b_n} + \log(1 + |a_n|) + 
k\, \E \log(1 + |U|) + (n-k)\, \E \log(1 + |V|) \\ 
 & \leq &  
\frac{1}{b_n} + \log(1 + |a_n|) + Cn. 
\ene 
Thus, by (4.7) and (4.6), with some constant $C$ 
$$ 
D(R_{k,n}||Z) \leq C \big(n + \log b_n + 1/b_n + \log(1 + |a_n|)\big). 
$$ 
It remains to apply this bound in (4.3)-(4.4), and then (4.1) yields 
$|\widetilde D_n - D_n| = o(b_1^n)$ with any $b_1 > b$. 
One may take $b_1 = \frac{1}{2}$, and thus Lemma 4.1 is proved.

 
\vskip5mm 
\section{{\bf Uniform Local Limit Theorem}} 
\setcounter{equation}{0} 
 
\vskip2mm 
Consider the normalized sums 
$Z_n = \frac{1}{b_n}\,(X_1 + \dots + X_n) - a_n$ (where $a_n \in \R, \ b_n > 0$), 
associated to independent identically distributed random variables $X_k$. 
 
\vskip5mm 
{\bf Proposition 5.1.} {\it Assume that $Z_n \Rightarrow Z$, 
where $Z$ has a $($continuous$)$ density $\psi$. If the random variables $Z_n$  
have absolutely continuous distributions for $n \geq n_0$ with densities, say  
$p_n$, then 
\be 
\sup_x \, |\widetilde p_n(x) - \psi(x)| \rightarrow 0 \qquad (n \rightarrow \infty). 
\en 
} 
 
Here $\widetilde p_n$ denote the modified densities of $p_n$, constructed  
in Section 3 for the case $n_0 = 1$. Necessarily, $Z$ has a stable distribution  
of some index  
$\alpha \in (0,2]$, and $b_n \sim n^{1/\alpha} B(n)$, where $B$ is a slowly 
varying function.  
 
Note that in this proposition 
it does not matter, whether $Z$ is extremal or not. 
 
\vskip5mm 
{\bf Proof.} Consider the characteristic functions 
$$ 
\widetilde f_n(t) = \int_{-\infty}^{+\infty} e^{itx} \widetilde p_n(x)\,dx,  
\qquad f(t) = \int_{-\infty}^{+\infty} e^{itx} \psi(x)\,dx, 
$$  
and express the densities via inverse Fourier transforms, while splitting  
the Fourier integral into the two regions, 
$$ 
\widetilde p_n(x) - \psi(x) = \frac{1}{2\pi} \int_{|t| \leq T_n} 
e^{-itx}\,(\widetilde f_n(t) - f(t))\,dt + 
\frac{1}{2\pi} \int_{|t| > T_n} e^{-itx}\,(\widetilde f_n(t) - f(t))\,dt 
$$ 
with given $t_0 > 0$ and $0 < T_n \leq t_0 b_n$.  
 
Let $f_n$ denote the characteristic functions of $Z_n$. By assumption, if  
$T_n \rightarrow +\infty$ sufficiently slowly, then  
$T_n \max_{|t| \leq T_n} |f_n(t) - f(t)| \rightarrow 0$. Hence, by (3.6), 
$$ 
T_n \max_{|t| \leq T_n} |\widetilde f_n(t) - f(t)| \rightarrow 0 \qquad  
(n \rightarrow \infty), 
$$ 
so that uniformly in all $x$ 
$$ 
\widetilde p_n(x) - \psi(x) = 
\frac{1}{2\pi} \int_{|t| > T_n} e^{-itx}\,(\widetilde f_n(t) - f(t))\,dt + o(1). 
$$ 
Moreover, since $f$ is integrable, 
\be 
\sup_x\, |\widetilde p_n(x) - \psi(x)| \leq  
\frac{1}{2\pi} \int_{|t| > T_n} |\widetilde f_n(t)|\,dt + o(1). 
\en 
 
Recall that the characteristic functions $\widetilde f_n$ are  
integrable as well. The integration in (5.2) should also be splitted into  
the two regions accordingly, and hence the integral itself will be bounded by 
$$ 
\int_{T_n < |t| < t_0 b_n} |\widetilde f_n(t)|\,dt +  
\int_{|t| > t_0 b_n} |\widetilde f_n(t)|\,dt. 
$$ 
But, by Lemma 3.2, the last integral tends to zero, as long as $b_n$ has  
at most polynomial growth. Using once more (3.6), we arrive at 
$$ 
\sup_x \, |\widetilde p_n(x) - \psi(x)| \leq 
\frac{1}{2\pi} \int_{T_n < |t| < t_0 b_n} |f_n(t)|\,dt + o(1). 
$$ 
 
It remains to apply the following bound derived in [I-L], p.133. 
There exist constants $c>0$ and $t_0 > 0$, such that 
$|f_n(t)| \leq e^{-c |t|^{\alpha/2}}$ for all $n \geq 1$ and $t$ in the interval  
$|t| < t_0 b_n$. This gives 
$$ 
\int_{T_n < |t| < t_0 b_n} |f_n(t)|\,dt \leq  
\int_{|t| < t_0 b_n} e^{-c |t|^{\alpha/2}}\,dt \rightarrow 0. 
$$ 
Thus, Proposition 5.1 is proved.

 
\vskip5mm 
\section{{\bf Proof of Theorem 1.1}} 
\setcounter{equation}{0} 
 
Let $X$ and $Z$ be random variables with densities $p$ and $\psi$, 
such that $\psi(x) = 0 \Rightarrow p(x) = 0$ a.e. The relative entropy 
\be 
D(X||Z) =  
\int_{-\infty}^{+\infty} \frac{p(x)}{\psi(x)} \log\frac{p(x)}{\psi(x)}\ \psi(x) dx 
\en 
is well-defined and may be bounded from above by applying an elementary  
inequality 
\be 
t\log t \leq (t-1) + C_\ep |t-1|^{1+\ep} \qquad (t \geq 0), 
\en 
where $C_\ep$ depends on $\ep \in (0,1]$, only. Namely, it immediately yields 
$$ 
D(X||Z) \leq C_\ep 
\int_{-\infty}^{+\infty} \frac{|p(x)-\psi(x)|^{1+\ep}}{\psi(x)^\ep}\, dx. 
$$ 
Moreover, letting $\Delta = \sup_x |p(x) - \psi(x)|$, 
\be 
D(X||Z) \leq C_\ep \Delta^\ep 
\int_{-\infty}^{+\infty} \frac{1}{\psi(x)^\ep}\,|p(x)-\psi(x)|\, dx. 
\en 
 
This is a general upper bound which may be used in the proof of Theorem 1.1 
in case of a non-normal stable density $\psi$ of index $0 < \alpha < 2$  
and using $X = \widetilde Z_n$ with modified densities $\widetilde p_n$. 
Indeed, by Proposition 5.1, 
$$ 
\Delta_n = \sup_x \, |\widetilde p_n(x) - \psi(x)| \rightarrow 0 \qquad  
(n \rightarrow \infty). 
$$ 
In addition, if $Z$ is non-extremal, $\psi$ admits a lower bound  
$\psi(x) \geq c\,(1+|x|)^{-(1+\alpha)}$ with some constant $c>0$, cf. (1.1).  
Hence, by (6.3), 
$$ 
D(\widetilde Z_n||Z) \leq C \Delta_n^\ep 
\int_{-\infty}^{+\infty} (1+|x|)^{\ep (1+\alpha)}\, 
|\widetilde p_n(x)-\psi(x)|\, dx, 
$$ 
where the constant depends on $\ep$ and $\psi$. But, for an arbitrary 
$\ep < \frac{\alpha}{1 + \alpha}$, so that $s = \ep (1+\alpha) < \alpha$, 
we see that the last integral does not exceed 
$$ 
\int (1+|x|)^s\,|\widetilde p_n(x)-p_n(x)|\, dx + 
\int (1+|x|)^s\,p_n(x)\, dx + 
\int (1+|x|)^s\,\psi(x)\, dx. 
$$ 
All these integrals are bounded by a constant, which follows from Lemma 3.1  
and the fact that $\sup_n \E\, |Z_n|^s < +\infty$ 
(which is due to the assumption $Z_n \Rightarrow Z$, cf. [I-L]). 
As a result, we have $D(\widetilde Z_n||Z) \rightarrow 0$, which yields  
the desired conclusion $D(Z_n||Z) \rightarrow 0$ in view of Lemma 4.1. 
 
In the normal case $(\alpha = 2)$, a similar argument with slight  
modifications may be applied as well. Without loss of generality, assume that 
$Z$ is standard normal, i.e., $\psi(x) = \frac{1}{\sqrt{2\pi}}\, e^{-x^2/2}$.  
Now we use (6.2) with $\ep = 1$ and $C_\ep = 1$.  
More precisely, splitting the integration in (6.1) into the two regions, 
we get 
\be 
D(X||Z) \leq 
\int_{|t| \leq T} \frac{|p(x)-\psi(x)|^2}{\psi(x)}\, dx + 
\int_{|t| > T} p(x) \log\frac{p(x)}{\psi(x)}\,dx 
\en 
with an arbitrary $T>0$. Furthermore, the first integral on the right-hand side 
is bounded by 
$$ 
e^{T^2/2} \sqrt{2\pi}\, \Delta \int_{|t| \leq T} |p(x)-\psi(x)|\, dx 
\leq 2\,e^{T^2/2} \sqrt{2\pi}\, \Delta. 
$$ 
If additionally $p(x) \leq M$, (6.4) leads to another general upper bound 
\be 
D(X||Z) \leq 2\,e^{T^2/2} \sqrt{2\pi}\, \Delta + 
\int_{|t| > T} \big(x^2 + \log(M\sqrt{2\pi}\big) p(x)\,dx. 
\en 
 
Here, again let $X = \widetilde Z_n$ and $T = T_n$. Then 
the above bound holds with $M = 1$ and all $n$ large enough. 
If $T_n \rightarrow +\infty$ sufficiently slow, we have 
$e^{T_n^2/2} \Delta_n \rightarrow 0$. On the other hand, 
\bee 
\int_{|t| > T_n} x^2\, \widetilde p_n(x)\,dx 
 & \leq &  
\int_{-\infty}^{+\infty} x^2\, |\widetilde p_n(x) - p_n(x)|\,dx \\ 
 & & + \, 
\int_{|t| > T_n} x^2\, p_n(x)\,dx + \int_{|t| > T_n} x^2\, \psi(x)\,dx. 
\ene 
Again, all the integrals tend to zero, in view of Lemma 3.1 and the uniform 
integrability of the sequence $Z_n^2$. Hence, by (6.5), 
$D(\widetilde Z_n||Z) \rightarrow 0$, which, by Lemma 4.1, proves  
Theorem 1.1 in the normal case.

\end{document}